\documentclass[12pt]{amsart}
\usepackage{latexsym}
\usepackage{amscd,amssymb,times,fullpage}
\newtheorem{thm}{Theorem}
\newtheorem{prop}{Proposition}
\newtheorem{lem}{Lemma}

\newtheorem{problem}{Problem}

\theoremstyle{remark}
\newtheorem{rem}{Remark}

\theoremstyle{definition}

\newcommand{\C}{\mathbb{ C}}

\newcommand{\bH}{\mathbb{ H}}
\newcommand{\R}{\mathbb{ R}}
\newcommand{\HH}{\mathcal{ H}}

\newcommand{\Z}{\mathbb{ Z}}

\title[Free circle actions with contractible orbits]{Free circle actions with contractible orbits on\\ 
symplectic manifolds}
\author{D.~Kotschick}
\address{Mathematisches Institut, Ludwig-Maximilians-Universit\"at M\"unchen,
Theresienstr.~39, 80333 M\"unchen, Germany}
\email{dieter@member.ams.org}
\date{March 28, 2005; MSC 2000: primary 57R17, 57R57, secondary 57M60, 57R91}

\begin{document}

\begin{abstract}
    We prove that closed symplectic four-manifolds do not admit any smooth free circle actions 
    with contractible orbits, without assuming that the actions preserve the symplectic forms. In higher
    dimensions such actions by symplectomorphisms do exist, and we give explicit examples
    based on the constructions of~\cite{FGM}.
\end{abstract}

\maketitle

\section{Introduction}

The following problem was raised by McDuff and Salamon in~\cite{MS}, p.~152:
\begin{problem}
Do there exist closed symplectic manifolds with symplectic free circle actions whose orbits are contractible?
\end{problem}
As remarked by McDuff and Salamon in {\it loc.~cit.}, it is easy to find examples where the orbits are null-homologous
but not null-homotopic. In fact some circle actions on the Kodaira--Thurston manifold have this property. It is also 
known, and is easy to see, that a free action can not be Hamiltonian. Even without assuming that the action is 
symplectic, we show in this paper that the answer to the Problem is negative in dimension four:
\begin{thm}\label{t:main}
There is no closed symplectic four-manifold that admits a smooth free circle action with contractible orbits.
\end{thm}
A proof is given in Section~\ref{s:main} below. Some variations in the argument are possible, depending on how 
much Seiberg--Witten theory and how much of the theory of three-manifolds one chooses to use. In an attempt 
to be as elementary as possible, our proof is almost purely topological, but it seems that Seiberg--Witten 
theory can only be avoided, if one does assume that the action is symplectic. In Section~\ref{s:FGM} we shall show 
how the arguments of Fernandez, Gray and Morgan~\cite{FGM} on the one hand give an elementary proof of 
Theorem~\ref{t:main} if one assumes the action to be symplectic, and, on the other hand, allow one to give a positive 
answer to the Problem in dimensions $\geq 6$:
\begin{thm}\label{t:higher}
In every even dimension $\geq 6$ there exist closed symplectic manifolds with symplectic free circle actions whose 
orbits are contractible.
\end{thm}
By the results of~\cite{FGM} all such manifolds are circle bundles over mapping tori. However, there are many different
choices one can make for the fiber of the mapping torus, and for its monodromy. We shall use non-trivial facts
about the geometry of $K3$ surfaces in order to produce concrete, explicit examples. 

In Section~\ref{s:final} we address a related question raised recently by Allday and Oprea~\cite{AO}. 
We show that $S^{3}\times S^{1}$ is the only compact complex surface which admits a free circle action with 
contractible orbits.

\section{Free circle actions with contractible orbits}\label{s:main}

In this section we prove Theorem~\ref{t:main}.
Let $X$ be a smooth manifold with a smooth free circle action with quotient $M$. 
\begin{lem}\label{lem1}
If the orbits of the circle action are contractible, then the projection $\pi\colon X\rightarrow M$ induces an 
isomorphism on fundamental groups, and an exact sequence
\begin{equation}\label{pi2}
1\longrightarrow\pi_2(X)\stackrel{\pi_*}{\longrightarrow}\pi_2(M)\stackrel{e}{\longrightarrow}\Z\longrightarrow 1 \ .
\end{equation}
The action of $\pi_1(M)$ on $\pi_2(M)$ preserves the subgroup $\pi_1(X)$.
\end{lem}
This follows immediately from the long exact sequence of homotopy groups for $\pi$. The connecting homomorphism 
$e$ is given by evaluation of the Euler class on spherical homology classes in $M$. 

\begin{prop}\label{split}
Let $X$ be an orientable smooth four-manifold with a smooth free circle action with quotient $M$. If the orbits are 
contractible, then $M$ is a connected sum $M_0\# M_1$, and $X$ is the fiber sum of two circle bundles 
$X_i\rightarrow M_i$, with $X_0=S^3\times S^1\rightarrow M_0=S^2\times S^1$ the product of the Hopf fibration 
with a trivial circle factor.
\end{prop}
\begin{proof}
Consider the Kneser--Milnor prime decomposition of $M$. The irreducible summands have vanishing $\pi_2$ by the 
sphere theorem, cf.~\cite{Hempel}, and the other summands are diffeomorphic to $S^2\times S^1$. By Lemma~\ref{lem1}
there are $k\geq 1$ summands diffeomorphic to $S^2\times S^1$, and applying Lemma~\ref{lem1} to the circle bundle over 
$k(S^2\times S^1)$, we see that there is a spherical homology class in $k(S^2\times S^1)$ on which the Euler 
class of $X\longrightarrow M$ evaluates as $+1$. This homology class is clearly primitive. Now the diffeomorphism group 
of $k(S^2\times S^1)$ surjects onto the automorphism group of the free group on $k$ generators (see p.~81 of~\cite{L}),
which in turn surjects onto $Gl(k,\Z)$. As $Gl(k,\Z)$ acts transitively on the primitive vectors in $\Z^k$, we see that every 
primitive homology class in $k(S^2\times S^1)$ is represented by an embedded $S^2$. Thus there is an embedded 
$2$-sphere $S\subset M$ on which the Euler class of the circle bundle evaluates as $+1$. Moreover, because the homology
class of $S$ is primitive, we can find an embedded $S^1\subset M$ intersecting $S$ transversely and precisely once.
The boundary of a regular neighbourhood of the union of $S$  and this $S^1$ is another $2$-sphere, splitting off 
$S^2\times S^1$ as a connected summand of $M$. The separating $2$-sphere in $M$ along which $S^2\times S^1$ 
splits off is null-homologous in $M$, so the Euler class evaluates trivially on it, and the circle bundle over it is trivial. 
Therefore in the total space $X$ of the fibration $\pi$, a copy of $S^3\times S^1$ splits off along $S^2\times S^1$, in such 
a way that the inclusion of $S^2\times S^1$  into $(S^3\times S^1)\setminus (D^3\times S^1)$ induces the zero map 
on fundamental groups. The fibration over this connected summand in the base is the Hopf fibration multiplied with 
the identity on the circle: 
\begin{equation}\label{Hopf}
S^3\times S^1\stackrel{\textrm{Hopf}\times\textrm{Id}}{\longrightarrow} S^2\times S^1 \ .
\end{equation} 
\end{proof}

As a converse to this proof, note that we could start with any fibration $X_1\rightarrow M_1$ and fiber sum it 
with~\eqref{Hopf} to obtain a circle bundle with homotopically trivial fibers. If $X_1\rightarrow M_1$ itself has the 
property that the fiber is homotopically trivial, then $X$ splits off two copies of~\eqref{Hopf} in a fiber sum. In this
situation we have:
\begin{lem}\label{diff}
The total space of the fiber sum of two copies of~\eqref{Hopf} is diffeomorphic to 
$(S^3\times S^1)\# (S^2\times S^2)\# (S^3\times S^1)$.
\end{lem}
\begin{proof}
This is implicit in the proof of Theorem~2 of~\cite{entropies}. Consider the standard action of $S^1$ on $S^2$
by rotation around the north-south axis. Its linearization at the two fixed points induces opposite orientations
of the tangent planes. Fixing an orientation of $S^2$, we assign a sign $=\pm 1$ to each fixed point according to
whether the linearization of the action induces the given orientation, or not. We can make an equivariant 
self-connected sum at the two fixed points, by identifying the boundaries of $S^1$-invariant neighbourhoods of 
the fixed points by an orientation-reversing $S^1$-equivariant diffeomorphism. This gives us $T^2$ with a standard 
$S^1$-action.

Now take the diagonal action on $S^2\times S^2$ corresponding to the above action on the factors. It has four 
fixed points $(p,q)$, where $p,q\in S^2$ are fixed points of the action on $S^2$. If $p$ and $q$ were assigned 
the same sign above, then the linearization of the diagonal action at $(p,q)$ induces the product orientation of 
$S^2\times S^2$ arising from the given orientation of $S^2$; if the signs were different, it induces the opposite 
orientation. Here too we can eliminate the fixed points by equivariant self-connected sum to obtain a free action 
on $(S^3\times S^1)\# (S^2\times S^2)\# (S^3\times S^1)$. The quotient of the diagonal action on $S^2\times S^2$
is $S^3$ with four marked points corresponding to the fixed points. Consider an embedded $2$-sphere $S\subset S^3$ 
which splits $S^3$ into two connected components, each containing a pair of marked points at which the linearized 
action induces opposite orientations. Then these points can be paired in the self-connected sum, and the preimage 
of $S$ in $(S^3\times S^1)\# (S^2\times S^2)\# (S^3\times S^1)$ is a copy of $S^2\times S^1$, splitting the fibration 
into the fiber sum of two copies of~\eqref{Hopf}.
\end{proof}

We can now prove Theorem~\ref{t:main}.

\begin{proof}[Proof of Theorem~\ref{t:main}]
Let $X$ be a smooth orientable four-manifold with a free circle action with contractible orbits.
Consider the splitting of $X$ into a fiber sum given by Proposition~\ref{split}. If $b_1(M_1)=0$, then 
$b_1(X)=1$, and the vanishing of the Euler characteristic of $X$ shows that $b_2(X)=0$, so that $X$ cannot be symplectic. 
If $b_1(M_1)>0$, then $\pi_1(M)=\pi_1(X)$ splits as a non-trivial free product in which
both free factors have Abelianizations of positive rank. Moreover, this splitting is realized by the connected sum 
decomposition of $M$, respectively the decomposition of $X$ as a fiber sum. We can therefore take four-fold covers 
of $X$ induced from four-fold covers of $M_1$, and then these covers split off four copies of~\eqref{Hopf}
in fiber sums. By Lemma~\ref{diff}, this means that the covers split off two copies of $S^2\times S^2$ as
connected summands. Therefore, the numerical Seiberg--Witten invariants of this cover are well-defined and vanish. 
If $X$ were symplectic then so would be any finite cover, leading to a contradiction with Taubes's non-vanishing 
result~\cite{T1,Bourbaki} for the Seiberg--Witten invariants of symplectic manifolds with $b_2^+>1$. 
\end{proof}

\section{Symplectic circle actions}\label{s:FGM}

Let $(X,\omega)$ be a symplectic manifold with a symplectic free circle action with quotient $M$. 
The $S^1$-invariant closed forms represent all cohomology classes, and nondegeneracy is an open condition in 
the space of closed invariant forms. Therefore an open neighbourhood of $[\omega]\in H^2(X;\R)$ can be represented 
by $S^1$-invariant symplectic forms. As rational points are dense, we may assume without loss of generality that the class 
$[\omega]$ is rational, and even integral.

Let $V$ be the vector field generating the action. Then $L_V\omega=0$, and the Cartan formula implies that 
$i_V\omega$ is a closed $1$-form. It is $S^1$-invariant and vanishes on the orbits, so that it is the pullback of a 
closed $1$-form $\alpha$ on $M$. This form is nowhere zero because $\omega$ is nondegenerate, and it represents 
an integral cohomology class because $\omega$ does. Thus integration of $\alpha$ defines a smooth fibration 
$M\longrightarrow S^1$.

Assume now that $M$ is three-dimensional. If the fiber of $M\longrightarrow S^1$ has positive genus, then $M$ is 
aspherical, and the orbits of the circle action on $X$ can not be contractible, see Lemma~\ref{lem1}. If the fiber of 
$M\longrightarrow S^1$ has genus zero, and  the orbits of the circle action are contractible, then $X$ fibers over
$S^1$ with fiber $S^3$. Therefore its second Betti number vanishes, contradicting the assumption that $X$ is 
symplectic. This proves Theorem~\ref{t:main} for symplectic circle actions.

Reversing this construction as in~\cite{FGM} allows us to prove Theorem~\ref{t:higher}. 

\begin{proof}[Proof of Theorem~\ref{t:higher}]
It suffices to construct a 
$6$-dimensional example. For this let $Y$ be the smooth oriented four-manifold underlying a complex $K3$ surface. 
Given that the intersection form of $Y$ is $Q=3H\oplus 2E_8$, we can choose non-zero primitive classes $x$ and 
$c\in H^2(Y;\Z)$ as follows.
Take $c$ to be a basis vector for the first copy of $H$ in the intersection form, and $x$ the sum of basis
vectors for the second copy of $H$. It is elementary to find an automorphism $f$ of $H\oplus H$ with $f(x)=x+c$ 
and $f(c)=c$. We then extend $f$ to all of $Q$ in such a way that it preserves an orientation of a maximal 
positive-definite subspace. Let $\HH_0$ be such a subspace containing $x$. Set $\HH_1=f(\HH_0)$, and let $\HH_t$ 
be a smoothly varying family of maximal positive definite subspaces with $t\in [0,1]$ interpolating between the given 
endpoints. We can choose this family so that $x+tc\in\HH_t$.

As the automorphism $f$ of the intersection form of $Y$ preserves the orientation of a maximal positive-definite 
subspace, one can find an orientation-preserving diffeomorphism $\Phi\colon Y\longrightarrow Y$ with $\Phi^*=f$ by
the result of Matumoto~\cite{M}, compare~\cite{D}. Let $M$ be the mapping torus of $\Phi$.

We construct differential forms on $M$ by gluing up forms on $Y\times [0,1]$ with $\Phi^*$. The resulting forms will be 
smooth if we only glue up forms which are constant near the ends of the interval. To arrange this fix once and for all a 
smooth monotonically increasing function $\psi\colon [0,1]\longrightarrow [0,1]$ with $\psi (0)=0$ and $\psi (1)=1$ which 
is constant near $0$ and near $1$. 

By a celebrated result of Yau~\cite{Y}, the maximal positive-definite subspace $\HH_{\psi (t)}$ for the intersection form 
determines a Ricci-flat K\"ahler--Einstein metric $g_{\psi (t)}$ for which $\HH_{\psi (t)}$ is the space of self-dual harmonic 
two-forms, cf.~\cite{BPV,D}.
This metric depends smoothly on $t$ because $\HH_{\psi (t)}$ does. As $x+\psi (t)c$ is contained in $\HH_{\psi (t)}$ by
construction, it is represented by a $g_{\psi (t)}$-self-dual harmonic form. The self-dual harmonic forms of  Ricci-flat 
K\"ahler--Einstein metrics are parallel, and therefore symplectic. Thus we have a family of symplectic two-forms
$\Omega_{\psi (t)}$ representing $x+\psi (t)c$. Now the diffeomorphism $\Phi$ is an isometry between $g_0$ and 
$g_1$, cf.~\cite{Y,BPV,D}. It maps $g_0$-self-dual harmonic two-forms to $g_1$-self-dual harmonic two-forms. Thus $f(x)=x+c$
and $f=\Phi^*$ imply $\Phi^*\Omega_0=\Omega_1$. Now we can glue up the path of two-forms $\Omega_{\psi (t)}$ 
on $Y$ to obtain a two-form $\Omega$ on the mapping torus $M$ which restricts 
to $Y\times\{ t\}$ as $\Omega_{\psi (t)}$.

The fact that $\Phi^*c=c$, implies that $c\in H^2(Y;\Z)$ lifts to a cohomology class on the mapping torus $M$ in the Wang
sequence
$$
\ldots\longrightarrow H^2(M;\Z)\longrightarrow H^2(Y;\Z)\stackrel{\Phi^*-1}{\longrightarrow} H^2(Y;\Z)\ldots \ .
$$
We choose such a lift and fix it once and for all. By abuse of notation we denote it by $c$. By Lemma~17 of~\cite{FGM} there
is a closed two-form $\gamma$ on $M$ representing the class $c$ and such that $d\Omega=\alpha\wedge\gamma$, where 
$\alpha$ is the pullback of the volume form on $S^1$ to $M$. 

Finally let $X\longrightarrow M$ be the circle bundle with Euler class $c$. Pick a connection $1$-form $\eta$ with $d\eta$ the 
pullback of $\gamma$. Let $\omega=\Omega+\alpha\wedge\eta$, where we have dropped the pullback from $M$ to $X$
in the notation. Then
$$
d\omega = d\Omega-\alpha\wedge d\eta = \alpha\wedge\gamma - \alpha\wedge\gamma = 0 \ .
$$
It is easy to see that $\omega$ is nondegenerate on $X$, and it is obviously $S^1$-invariant. Thus it is an invariant
symplectic form.

Now $\pi_2(M)=\pi_2(Y)$, and all homology classes in $H_2(Y;\Z)$ are spherical. As we chose $c$ to be primitive, there is a 
homology class on which $c$ evaluates as $+1$. As this class is spherical, we conclude that the evaluation of the Euler class 
$\pi_2(M)\stackrel{e}{\longrightarrow}\pi_1(S^1)$ in the exact homotopy sequence of $S^1\longrightarrow X\longrightarrow M$
is surjective. Thus the fiber  of this circle bundle is contractible in the total space. This proves Theorem~\ref{t:higher} in
dimension $6$, and the higher-dimensional case follows by taking products.
\end{proof}

\section{Circle actions on complex surfaces}\label{s:final}

Recently Allday and Oprea~\cite{AO} gave an example of a four-manifold with a free circle action with 
contractible  orbits which is $c$-symplectic, meaning that there is a degree $2$ cohomology class $c$ with $c^2>0$. 
Their example is constructed by doubling the complement of a regular neighbourhood of a fiber in~\eqref{Hopf}. 
We can think of this manifold as the fiber sum of two copies of~\eqref{Hopf}, if we reverse the orientation on one
of the copies before performing the fiber sum. Therefore Lemma~\ref{diff} shows that the example of Allday and 
Oprea~\cite{AO} is diffeomorphic to $(S^3\times S^1)\# (S^2\times S^2)\# (S^3\times S^1)$. Such examples
had also appeared in~\cite{entropies}, with a very different motivation. The proof of Theorem~2 in~\cite{entropies}
shows that one can change the smooth structure on a manifold with a free circle action so that for the new smooth 
structure there is no fixed-point-free circle action.

We found Theorem~\ref{t:main} generalizing the observation that the example of Allday and 
Oprea~\cite{AO} can not have a symplectic structure. In their paper, these  authors also asked whether their
example has a complex or an almost complex structure. Now an orientable four-manifold with a free circle action 
is parallelizable, because the quotient three-manifold is. Therefore, such manifolds are always almost complex.
But, contractibility of the orbits rules out complex structures in almost all cases:
\begin{thm}\label{t:cx}
Let $X$ be a compact complex surface with a free circle action with contractible orbits. Then $X$ is a primary
Hopf surface diffeomorphic to $S^{3}\times S^{1}$.
\end{thm}
\begin{proof}
Consider the fiber sum decomposition of $X$ into 
$X_0=S^3\times S^1\stackrel{\textrm{Hopf}\times\textrm{Id}}{\longrightarrow} S^2\times S^1=M_0$ and 
$X_1\longrightarrow M_1$ given by Proposition~\ref{split}. If $M_{1}$ is simply connected, then 
$\pi_{1}(X)=\Z$ and $b_{2}(X)=0$. A theorem of Kodaira then implies that $X$ is a primary Hopf surface, 
i.~e.~it is diffeomorphic to $S^{3}\times S^{1}$, cf.~\cite{N}.
If $M_1$ is not simply connected, then we have a decomposition of $\pi_1(X)$ into a non-trivial free product
$\Z * \Gamma$, implying that $\pi_1(X)$ has infinitely many ends. This is impossible by the following 
Proposition, which completes the proof.
\end{proof}
\begin{prop}
The fundamental group of a compact complex surface with vanishing Euler characteristic has one or two ends.
\end{prop}
\begin{proof}
Let $X$ be a compact complex surface with vanishing Euler number.
If $X$ is K\"ahlerian\footnote{By Theorem~\ref{t:main}, this case is not relevant for Theorem~\ref{t:cx}.}, 
then its fundamental group has at most one end, cf.~\cite{ABCKT}.

If $X$ is not K\"ahlerian, then by the Kodaira classification~\cite{BPV,ABCKT}, its first Betti number is odd, and 
$X$ is either elliptic, or of class VII. If $X$ is elliptic with odd $b_1(X)$ and vanishing Euler number, then its 
universal cover is diffeomorphic to $\R^4$ or to $S^3\times\R$ by a result of Kodaira, cf.~\cite{ABCKT,N}. 
Thus $\pi_1(X)$ has one or two ends.

If $X$ is of class VII, then $b_1(X)=1$. 
The vanishing of the Euler number implies $b_2(X)=0$. Now if $X$ admits non-constant meromorphic functions, 
then it is elliptic, by another result of Kodaira, cf.~\cite{N}, and this case has been dealt with already. We may 
therefore assume that $X$ has no non-constant meromorphic functions. Then, if $X$ contains a holomorphic curve, 
it is a Hopf surface with universal cover $\C^{2}\setminus\{ 0\}$, cf.~\cite{N}. If it contains no holomorphic curve, 
then it is an Inoue surface with universal cover $\bH\times\C$, cf.~\cite{K,N}. Again we conclude that $\pi_{1}(X)$ 
has one or two ends.
\end{proof}
\begin{rem}
The assumption that $X$ has vanishing Euler number was only used to bypass the lack of a classification of 
class VII surfaces with positive Euler number. In the K\"ahler case the value of the Euler number is irrelevant.
\end{rem}

\medskip
\noindent
{\bf Acknowledgements:} I am very grateful to R.~Fintushel and S.~Kerckhoff for useful discussions, 
and to an anonymous referee who requested clarification of the proof of Proposition~\ref{split}.
I would also like to thank Y.~Eliashberg for his hospitality at Stanford University during the preparation of this paper.

\bigskip

\bibliographystyle{amsplain}

\bigskip

\end{document}